\documentclass[preprint,12pt]{elsarticle}

\usepackage{amssymb}
\journal{arXiv}
\begin{document}

\begin{frontmatter}

\title{Two-sample extended empirical likelihood for estimating equations}
\author{Min Tsao and Fan Wu}
\address{Department of Mathematics and Statistics, University of Victoria, Victoria, British Columbia, Canada V8W 3R4}

\begin{abstract}

We propose a two-sample extended empirical likelihood for inference on the difference between two $p$-dimensional parameters defined by estimating equations. The standard two-sample empirical likelihood for the difference is Bartlett correctable but its domain is a bounded subset of the parameter space. We expand its domain through a composite similarity transformation to derive the two-sample extended empirical likelihood which is defined on the full parameter space. The extended empirical likelihood has the same asymptotic distribution as the standard one and can also achieve the second-order accuracy of the Bartlett correction. 
%A simple first-order version of this extended empirical likelihood is easy to compute and we recommend it for practical applications. 
We include two applications to illustrate the use of two-sample empirical likelihood methods and to demonstrate the superior coverage accuracy of the extended empirical likelihood confidence regions.
\\[0.2in]
\end{abstract}


\end{frontmatter}

\section{Introduction}  
A two-sample problem is concerned with making inference for the difference between the corresponding parameters of two populations/models with two independent samples. The difference between two population means is a special case that has been extensively studied; when the sample sizes are not large and underlying distributions are normal, methods for the Behren-Fisher problem or a two-sample $t$ method can be used; when the sample sizes are large, non-parametric $z$ based procedures can be used. Recently, the empirical likelihood method (Owen, 2001) has been successfully applied to this special case. See Jing (1995), Liu et al. (2008), Liu and Yu (2010), Wu and Yan (2012) and Wu and Tsao (2013). These empirical likelihood methods complement existing methods as they do not require strong conditions and are more accurate than normal approximation based methods when the underlying distributions are skewed. In particular, the extended two-sample empirical likelihood for the difference between two $p$-dimensional means (Wu and Tsao, 2013) is defined on the whole of $\mathbb{R}^p$ and is more accurate than other empirical likelihood methods.

In this paper, we study empirical likelihood methods for the general two-sample problem concerning the difference between two $p$-dimensional parameters defined by general estimating equations. The main contribution of this paper is a new extended empirical likelihood for such a difference. 
%We work in a non-parametric or semi-parametric sitting where the underlying distributions are not specified or not fully specified. 
The empirical likelihood method was introduced by Owen (1988, 1990). It has since been applied to many problems in statistics; see Owen (2001) and references therein. In particular, Qin and Lawless (1994) showed that the empirical likelihood is effective for inference on parameters defined by estimating equations. DiCiccio et al. (1991) and Chen and Cui (2007) proved that the empirical likelihood for estimating equations is Bartlett correctable; the Bartlett corrected empirical likelihood enjoys the second-order accuracy. Although there have been relatively few publications that apply empirical likelihood to the general two-sample problem, it is well-suited for this problem as the formulation of the one-sample empirical likelihood for estimating equations can be readily extended to handle the two-sample case; see, {\em e.g.,} Jing (1995), Qin and Zhao (2000), Liu et al. (2008), Liu and Yu (2010), Wu and Yan (2012) and Zi et al. (2012). In particular, Qin and Zhao (2000) studied the standard two-sample empirical likelihood for the univariate version ($p=1$) of the problem, and Zi et al. (2012) considered the special case where the parameters are the coefficient-vectors of two linear models. 

In Section 2, we study the standard two-sample empirical likelihood for estimating equations in the general multi-dimensional sitting where $p\geq 1$. Like its one-sample counterpart, this two-sample empirical likelihood also has an asymptotic chi-square distribution and is Bartlett correctable.  Adopting the terminology in Tsao and Wu (2013), we refer to this standard two-sample empirical likelihood as the two-sample original empirical likelihood (OEL) for estimating equations. The OEL suffers from a mismatch problem (Tsao and Wu, 2013) in that it is only defined on a part of the parameter space. This problem affects the coverage accuracy of the OEL based confidence regions. To overcome this, in Section 3 we introduce a two-sample extended empirical likelihood (EEL) that is defined on the whole parameter space. The EEL is obtained by expanding the domain of the OEL to the full parameter space through a composite similarity mapping. We show that the EEL has the same asymptotic chi-square distribution as the OEL and that it can also achieve the second-order accuracy of the Bartlett correction. In Section 4, we discuss two applications of the two-sample OEL and EEL. The first application is concerned with the inference for the difference between two Gini indices, and the second application is concerned with that between coefficient vectors of two regression models. We also make use of these applications to compare the numerical accuracy of the OEL and EEL confidence regions and to illustrate the superior accuracy of the EEL.

Proofs of theoretical results on two-sample OEL and EEL are all relegated to the Appendix. Note that some of these results can be proved by slightly modifying the proofs of similar results for other empirical likelihoods in the literature. For brevity, we will not include detailed proofs for such results in the Appendix but will give relevant references which contain proofs of similar results. 

%2222222222222222222222222222222222222222222222222222222222222222222222222222222222
\section{Two-sample Original empirical likelihood (OEL) for estimating equations}

We first describe the general two-sample problem for estimating equations as follows.
Let $\{X_1, \dots, X_m\}$ and $\{Y_1, \dots, Y_n\}$ be independent copies of random vectors $X\in \mathbb{R}^{d}$ and $Y\in \mathbb{R}^{d}$ with parameters $\theta_x\in \mathbb{R}^p$ and $\theta_y\in \mathbb{R}^p$, respectively. Let $g(X,\theta_x)$ and $g(Y,\theta_y)$ be two $q$-dimensional estimating functions for $\theta_{x_0}$ and $\theta_{y_0}$ satisfying $E\{g(X,\theta_{x_0})\}=0$ and $E\{g(Y,\theta_{y_0})\}=0$, respectively. The unknown parameter of interest is the difference $\pi_0=\theta_{y_0}-\theta_{x_0}\in \mathbb{R}^{p}$ and the parameter space is the entire $\mathbb{R}^{p}$. A more general version of this problem allows the estimating function for $\theta_{x_0}$ to be different from that for $\theta_{y_0}$. For simplicity, we consider only the common case where the two estimating functions are the same.

We now generalize the one-sample OEL for estimating equations (Qin and Lawless, 1994) to obtain a two-sample OEL for $\pi_0$ and study its asymptotic properties.
We will need the following four conditions on $g(X,\theta_x)$ and $g(Y,\theta_y)$.  \vspace{0.1in}

{\em Condition 1.} $E\{g(X,\theta_{x_0})\}=0$ and $E\{g(Y,\theta_{y_0})\}=0$, and ${\rm var}\{g(X,\theta_{x_0})\}\in \mathbb{R}^{q\times q}$ and  ${\rm var}\{g(Y,\theta_{y_0})\}\in \mathbb{R}^{q\times q}$ are both positive definite.

{\em Condition 2.} $\partial g(X,\theta_x)/\partial \theta_x$ and $\partial g^2(X,\theta_x)/\partial \theta_x \partial \theta_x^T$ are continuous in $\theta_x$, and for $\theta_x$ in a neighbourhood of $\theta_{x_0}$ they are each bounded in norm by an integrable function of $X$. 

{\em Condition 3.} $\partial g(Y,\theta_y)/\partial \theta_y$ and $\partial g^2(Y,\theta_y)/\partial \theta_y \partial \theta_y^T$ are continuous in $\theta_y$, and for $\theta_y$ in a neighbourhood of  $\theta_{y_0}$ they are each bounded in norm by an integrable function of $Y$. 

{\em Condition 4.} $\lim \sup_{\|t\|\rightarrow \infty}|E[\exp\{it^Tg(X,\theta_x)\}]| <1$ and $E\|g(X,\theta_x)\|^{15}<+\infty$; $\lim \sup_{\|t\|\rightarrow \infty}|E[\exp\{it^Tg(Y,\theta_y)\}]| <1$ and $E\|g(Y,\theta_y)\|^{15}<+\infty$.  

\vspace{0.1in}

Denote by $\bar{p}=(p_1,...,p_m)$ and $\bar{q}=(q_1,...,q_n)$ two probability vectors satisfying $p_i\geq 0$, $q_j\geq 0$, $\sum_{i=1}^m p_i=1$ and $\sum_{i=1}^n q_j=1$. Let $\theta_y$ and $\theta_x$ be points in $\mathbb{R}^p$ and denote by $\theta_y(\bar{q})$ and $\theta_x(\bar{p})$ values that satisfy
\[\sum_{i=1}^mp_ig(X_i,\theta_x(\bar{p}))=0, \hspace{0.2in} \sum_{j=1}^np_jg(Y_j,\theta_y(\bar{q}))=0.\]

Let $\pi = \theta_y-\theta_x \in \mathbb{R}^p$ and let $\pi(\bar{p},\bar{q})=\theta_y(\bar{q})-\theta_x(\bar{p})$. Then, the two-sample OEL for a possible value of the difference $\pi$, $L(\pi)$, is defined as
\begin{equation}
L(\pi) = \sup_{(\bar{p},\bar{q}):\pi(\bar{p},\bar{q})=\pi} \left(\prod^{m}_{i=1}p_i \right) \left(\prod^{n}_{j=1}q_j \right), \label{oel}
\end{equation}
which is the maximum of the product of the one-sample OEL for $\theta_y$ and the one-sample OEL for $\theta_x$ taken over all pairs $(\theta_x, \theta_y)$ that satisfies $\pi = \theta_y-\theta_x$.
The corresponding two-sample empirical log-likelihood ratio for $\pi$ is thus
\begin{equation}
l(\pi) = -2\sup_{(\bar{p},\bar{q}):\pi(\bar{p},\bar{q})=\pi} \left\{ \sum^{m}_{i=1}\log(mp_i)+ \sum^{n}_{j=1}\log(nq_j) \right\}. \label{oel.log}
\end{equation}
For convenience, we will also use OEL for the original empirical log-likelihood ratio. We will write ``OEL $L(\pi)$'' and ``OEL $l(\pi)$'' for the original empirical likelihood ratio and log-likelihood ratio, respectively. 

Define the domain of $L(\pi)$, $D_n$, as
\begin{center}
$D_{n}=\{\pi\in \mathbb{R}^p: \mbox{there exist $\theta_y(\bar{q})$ and $\theta_x(\bar{p})$ such that $\pi=\theta_y(\bar{q})-\theta_x(\bar{p}). $}$\},
\end{center}
and define the domain of $l(\pi)$, $\Pi_n$, as
\begin{center}
$\Pi_n=\{\pi: \pi\in D_{n}\mbox{ and } l(\pi) < +\infty \}.$
\end{center}

Let $N=m+n$, $f_m=N/m$ and $f_n=N/n$. Without loss of generality, we assume that $m\geq n>q$. We also assume that $m/n=O(1)$ so that $O(n^{-1})$, $O(m^{-1})$ and $O(N^{-1})$, for example, are all interchangeable. By the method of Lagrangian multipliers, we have
\begin{equation}
l(\pi_0) = 2\left[\sum^{n}_{j=1} \log \{1+f_n (\lambda^*)^{T}g(Y_j, \theta^*_{y})\}+  \sum^{m}_{i=1} \log \{1-f_m (\lambda^*)^{T}g(X_i, \theta^*_{x})\} \right] \label{oel.lag.mul}
\end{equation}
where $(\lambda^*,\theta^*_y,\theta^*_x)$ is the solution of the following non-linear system
\begin{eqnarray}\label{oel.xy}
   && \sum^{n}_{j=1} \frac{g(Y_j, \theta_{y})}{1+f_n \lambda^{T}g(Y_j, \theta_{y})}=0, \\  \nonumber%\hspace{0.2in} %\mbox{and} \hspace{0.2in}
   &&\sum^{m}_{i=1} \frac{g(X_i, \theta_{x})}{1-f_m \lambda^{T}g(X_i, \theta_{x})}=0,  \\   
   &&\pi_0=\theta_y-\theta_x.   \nonumber 
\end{eqnarray}
Hence, we may write $l(\pi_0)=l(\lambda^*,\theta^*_y,\theta^*_x)$. 
The following theorem gives the asymptotic distribution of $l(\pi_0)$.

\vspace{0.1in}

\noindent {\bf Theorem 1.} \emph{
Under Conditions 1, 2, 3 and 4, the two-sample OEL $l(\pi_0)$ defined by (\ref{oel.lag.mul}) satisfies
\begin{equation}
l(\pi_0) \stackrel{D}{\longrightarrow}\chi^2_q \hspace{0.2in} \mbox{as} \hspace{0.1in} n\rightarrow +\infty.  \label{oel.chi}
\end{equation}
}

\vspace{0.1in}

By Theorem 1, the 100($1-\alpha$)\% two-sample OEL confidence region for $\pi_0$ is
\begin{equation}
{\mathcal C}_{1-\alpha}=\{\pi: \pi \in \mathbb{R}^{p} \mbox{ and } l(\pi) \leq c_\alpha \}   \label{oel.region}
\end{equation}
where $c_\alpha$ is ($1-\alpha$)th quantile of the $\chi^2_q$ distribution. The coverage error of $\mathcal{C}_{1-\alpha}$ is $O(n^{-1})$, that is
\begin{equation}
P(\pi_0\in {\mathcal C}_{1-\alpha})= P\{l(\pi_0) \leq c_\alpha\}=1-\alpha + O(n^{-1}). \label{oel.coverage}
\end{equation}

%The proof of one-sample empirical likelihood for estimating equation was given by Qin and Lawless (1994) which is similar to Theorem 1 for two-sample OEL, hence omitted here. 
Theorem 1 is the standard first-order result for an OEL. The error rate of $O(n^{-1})$ in (\ref{oel.coverage}) follows from an argument in DiCiccio et al. (1991) for that of the one-sample empirical likelihood. See also Hall and La Scala (1990).

%Since Theorem 1 is implied by Theorem 2 of this section which gives the stronger second-order result, for brevity the proof of Theorem 1 is omitted. 

DiCiccio et al. (1991) and Chen and Cui (2007) showed that the one-sample OEL for estimating equations is Bartlett correctable; the Bartlett correction reduces the coverage error of the empirical likelihood confidence region to $O(n^{-2})$. Theorem 2 shows that the two-sample OEL for estimating equations (\ref{oel.lag.mul}) is also Bartlett correctable. The key result for proving Theorem 2 is Lemma 1 below. In order to 
present Lemma 1, we need to first introduce some new notations. 

Denote by $\theta^k_{y}$, $\theta^k_{x}$ and $\zeta^k$ approximations of $\theta^*_y,\theta^*_x$ and $\lambda^*$, respectively. For brevity, the analytic expressions of $\theta^k_{y}$, $\theta^k_{x}$ and $\zeta^k$ will be given later in the Appendix. For these three notations, we note that the $k$ in say $\theta^k_{y}$ is a superscript (not to the power of $k$) which indicates the order of the approximation is $O(n^{-(k+1)/2})$, {\em i.e.}, $\theta^k_{y}=\theta_{y}+O(n^{-(k+1)/2})$.  Let $V_1=f_n {\rm var}\{g(Y,\theta_{y})\}$, $V_2=f_m {\rm var}\{g(X,\theta_x)\}$, $V=V_1+V_2$ and $W=V_1V^{-1}V_2$. Further, define
\begin{eqnarray}
&&z_{j0}=V^{-1/2}g(y_j, \theta^0_{y}),  \hspace{0.1in}  z_{i0}=V^{-1/2}g(x_i, \theta^0_x) \hspace{0.1in} z_{j1}=V^{-1/2}g(y_j, \theta^1_{y}), \\ \nonumber 
&&z_{i1}=V^{-1/2}g(x_i, \theta^1_{x}), \hspace{0.1in} 
 s^{t_1t_2\dots t_l}=f_n^{l-1}E(z^{t_1}_{j0}z^{t_2}_{j0}\dots z^{t_l}_{j0})+(-1)^{l}f_m^{l-1}E(z^{t_1}_{i0}z^{t_2}_{i0}\dots z^{t_l}_{i0}),\\ \nonumber  
&& S^{t_1t_2\dots t_l}=\frac{f_n^{l-1}}{n}\sum^{n}_{j=1}(z^{t_1}_{j0}z^{t_2}_{j0}\dots z^{t_l}_{j0})+\frac{(-1)^{l}f_m^{l-1}}{m} \sum^{m}_{i=1}(z^{t_1}_{i0}z^{t_2}_{i0}\dots z^{t_l}_{i0})- s^{t_1t_2\dots t_l}, \\ \nonumber 
&& S^{t_1t_2\dots t_l}_1=\frac{f_n^{l-1}}{n}\sum^{n}_{j=1}(z^{t_1}_{j1}z^{t_2}_{j1}\dots z^{t_l}_{j1})+\frac{(-1)^{l}f_m^{l-1}}{m} \sum^{m}_{i=1}(z^{t_1}_{i1}z^{t_2}_{i1}\dots z^{t_l}_{i1})- s^{t_1t_2\dots t_l}, \nonumber %\hspace{0.1in} 
\end{eqnarray}
and 
\begin{eqnarray}
\nonumber  \Delta_1&=&S^{\tau}S^{\tau}-S^{\tau\upsilon}S^{\tau}S^{\upsilon}+\frac{2}{3}s^{\tau\alpha\beta}S^{\tau}S^{\alpha}S^{\beta}+S^{\tau\upsilon}S^{\upsilon\omega}S^{\tau}S^{\omega}+\frac{2}{3}S^{\tau\alpha\beta}S^{\tau}S^{\alpha}S^{\beta},\\  \nonumber       &&-2s^{\tau\upsilon\omega}S^{\tau\alpha}S^{\upsilon}S^{\alpha}S^{\upsilon}+s^{\tau\upsilon\omega}s^{\tau\alpha\beta}S^{\upsilon}S^{\omega}S^{\alpha}S^{\beta}-\frac{1}{2}s^{\tau\upsilon\omega\alpha}S^{\tau}S^{\upsilon}S^{\omega}S^{\alpha}, \\ \nonumber
\Delta_2 &=&(S^{\tau\upsilon}-S_1^{\tau\upsilon})S^{\tau}S^{\upsilon}+\frac{2}{3}(S_1^{\tau\alpha\beta}-S^{\tau\alpha\beta})S^{\tau}S^{\alpha}S^{\beta}. \nonumber
\end{eqnarray}
where we have used the common summation convention that if an index appears more than once in an expression, summation over the index is understood.

\vspace{0.1in}

\noindent {\bf Lemma 1.} \emph{
With above notations and under condition 1, 2, 3, and 4, we have
\begin{equation}
\frac{l(\pi_0)}{N} = \Delta_1+\Delta_2+O_p(n^{-5/2}). \label{oel.bart}
\end{equation}
}

\vspace{0.1in}

To see the connection between expansion (\ref{oel.bart}) and that of other high-order expansions of empirical log-likelihood ratios, we note that the $\Delta_1$ term in (\ref{oel.bart}) is similar to the expansion of the one-sample empirical log-likelihood ratio at the true parameter value given by DiCiccio et al. (1991) and Chen and Cui (2007). In the present case, the expansion at the true difference $\pi_0$ depends on the true parameter values $\theta_{x0}$ and $\theta_{y0}$, both of which need to be estimated. The use of the estimated values of these parameters resulted in the extra term $\Delta_2$ in expansion (\ref{oel.bart}). See also a similar $\Delta_2$ term in the expansion of the two-sample empirical log-likelihood ratio for the mean in Liu et al. (2008).

We now use Lemma 1 to derive the two-sample Bartlett corrected empirical likelihood confidence region for the difference between two parameters defined by estimating equations. Let $\eta$ be the Bartlett correction factor where
\begin{equation}
\eta=-\frac{1}{3d}s^{\tau\upsilon\omega}s^{\tau\alpha\beta}+\frac{1}{2d}s^{\tau\tau\alpha\alpha}+\frac{f_mf_n}{d}tr(V^{-1/2}WV^{-1/2}). \label{oel.bartfactor}
\end{equation} 
The derivation of (\ref{oel.bartfactor}) is similar to that for the Bartlett correction factor for the difference between two means in Liu et al. (2008), which involves taking expectations of $\Delta_1$ and $\Delta_2$ and omitting terms of order $O(n^{-1})$. With $\eta$, the two-sample Bartlett corrected empirical log-likelihood ratio (BEL) is given by 
$$l_B(\pi)=l(\pi)(1-\eta N^{-1}).$$
It follows that the two-sample BEL confidence region ${\mathcal C}_{1-\alpha}'$ for $\pi_0$ is 
\begin{equation}
{\mathcal C}_{1-\alpha}'= \{\pi: \pi \in \mathbb{R}^{p} \mbox{ and } l_B(\pi) \leq c \}. \label{bel.region}
\end{equation} 
Theorem 2 below shows the coverage error of ${\mathcal C}_{1-\alpha}'$ is $O(n^{-2})$.

\vspace{0.1in}

\noindent {\bf Theorem 2.} \emph{
Under Conditions 1, 2, 3 and 4, for any $c>0$ the Bartlett corrected two-sample empirical likelihood confidence region satisfies
\begin{equation}
P(\pi_0 \in {\mathcal C}'_{1-\alpha})=P[l(\pi_0)\{1-\eta N^{-1}\}\leq c]=P(\chi_d^2\leq c) +O(n^{-2}).\label{oel.bart.00}
% P(\pi_0 \in {\mathcal C}'_{1-\alpha})=P[l(\pi_0)\{1-\eta N^{-1}\}\leq c]=P(\chi_d^2\leq c) +O(n^{-2}).\label{oel.bart.00}
\end{equation}
}

\vspace{0.1in}

A stronger result due to DiCiccio et al. (1991) is that
\begin{equation}
P[l(\pi_0)\{1-\eta N^{-1}+O_p(n^{-3/2})\}\leq c]=P(\chi_d^2\leq c) +O(n^{-2}).\label{oel.bart.01}
% P(\pi_0 \in {\mathcal C}'_{1-\alpha})=P[l(\pi_0)\{1-\eta N^{-1}\}\leq c]=P(\chi_d^2\leq c) +O(n^{-2}).\label{oel.bart.00}
\end{equation}

%Also, the proof of (\ref{oel.bart.00}) depends critically on Lemma 1. See the Appendix.

%With Lemma \ref{lemma1}, the proof of Theorem 2 follows steps that led to formula (2.5) of DiCiccio et al. (1991) outlined in that paper. Thus for brevity, the proof of Theorem 2 is omitted.

The Bartlett correction factor $\eta$ in (\ref{bel.region}) and (\ref{oel.bart.00}) depends on the moments of $g(X;\theta_{x0})$ and $g(Y;\theta_{y0})$ which are not available in empirical likelihood applications. Fortunately, by (\ref{oel.bart.01}) we can use a $\surd{n}$-consistent estimator $\hat{\eta}$ in place of the $\eta$ in (\ref{oel.bart.00}) without affecting the $O(n^{-2})$ error term in (\ref{oel.bart.00}). In real applications of the Bartlett correction, the $\surd{n}$-consistent estimator $\hat{\eta}$ is usually used instead of the exact $\eta$; see for example, Chen and Cui (2007) and Liu and Chen (2010). For two-sample BEL for the difference between two means, Liu et al. (2008) gave a moment estimator $\hat{\eta}$ for $\eta$. Liu and Yu (2010) reported that $\hat{\eta}$ tends to underestimate $\eta$ and proposed a less biased estimator $\tilde{\eta}$ for $\eta$. This less biased $\tilde{\eta}$ is also applicable to our two-sample BEL for $\pi_0$, and we will use this $\tilde{\eta}$ for our simulation studies.

%3333333333333333333333333333333333333333333333333333333333333333333333333333333333333333333333333
\section{Two-sample extended empirical likelihood (EEL) for estimating equations}

\subsection{Composite similarity mapping}
Like the one-sample OEL for estimating equations, the two-sample OEL $l(\pi)$ also suffers from the mismatch problem between its domain $\Pi_n$ and the parameter space since the parameter space is $\mathbb{R}{^p}$ but $\Pi_n\subset \mathbb{R}{^p}$. The mismatch problem is a main contributor to the undercoverage problem of the OEL confidence regions. To solve this problem, we now expand $\Pi_n$ to match the parameter space $\mathbb{R}{^p}$ through a {composite similarity mapping} (Tsao and Wu, 2013). This leads to an EEL defined on $\mathbb{R}{^p}$ and hence is free from the mismatch problem.

Denote by $\tilde{\theta}_{x}$ and $\tilde{\theta}_{y}$ the $\surd{n}$-consistent maximum empirical likelihood estimators (MELEs) for $\theta_{x0}$ and $\theta_{y0}$, respectively. Then, it is not difficult to show that the MELE of $\pi_0$ is $\tilde{\pi}=\tilde{\theta}_{y}-\tilde{\theta}_{x}$ which is $\surd{n}$-consistent for $\pi_0$. We define the composite similarity mapping $h^C_N: \Pi_n \rightarrow \mathbb{R}^p$ as
\begin{equation}
h_N^C(\pi)=\tilde{\pi}+\gamma\{N,l(\pi)\}(\pi-\tilde{\pi}) \hspace{0.2in} \mbox{for $\pi \in \Pi_n$}, \label{h.function}
\end{equation}
where function $\gamma\{N,l(\pi)\}$ is the expansion factor given by the following expression which depends continuously on $\pi$ 
\begin{equation}
\gamma\{N,l(\pi)\}=1+\frac{l(\pi)}{2N}. \label{gamma.f}
\end{equation}
To see how $h^C_N$ maps $\Pi_n$ onto $\mathbb{R}^p$, define the level-$\tau$ OEL contour as
\begin{equation}
c(\tau)=\{\pi: \pi \in \Pi_n \hspace{0.07in} \mbox{and} \hspace{0.07in} l(\pi)=\tau\}, \label{3.10}
\end{equation}
where $\tau\geq \tilde{\tau}=l(\tilde{\pi})\geq 0$. For the just-determined case, the one-sample OEL's satisfy $l(\tilde{\theta}_{x})=1$ and $l(\tilde{\theta}_{y})=1$. Thus,
$L(\tilde{\pi})=1$ and $\tilde{\tau}=l(\tilde{\pi})=0$. The contours form a partition of the domain $\Pi_n$; that is,
$c(\tau_1) \cap c(\tau_2)=\varnothing$ for any $\tau_1 \neq \tau_2$ and
\begin{equation}
\Pi_n = \bigcup_{\tau\in[\tilde{\tau},+\infty)} c(\tau). \label{3.11}
\end{equation}
In addition to conditions 1 to 4 above, we now introduce a new condition.

{\em Condition 5.} Each contour $c(\tau)$ is the boundary of a connected region in $\mathbb{R}{^p}$, and the contours are nested in that if $\tau_1<\tau_2$, then $c(\tau_1)$ is contained in the interior of the region defined by $c(\tau_2)$.

Under Condition 5 and in view of (\ref{3.11}), the MELE $c(\tilde{\tau})=\{\tilde{\pi}\}$ may be regarded as the centre of domain $\Pi_n$. It follows that the value of $\tau$ measures the outwardness of a $c(\tau)$ with respect to the centre; the larger the $\tau$ value, the more outward $c(\tau)$ is. 
The following theorem gives three key properties of $h^C_N$. \vspace{0.1in}

\vspace{0.1in}

\noindent {\bf Theorem 3.} \emph{
Under conditions 1, 2 and 3, mapping $h^C_N$ defined by (\ref{h.function}) and (\ref{gamma.f}) satisfies 
(i) it has a unique fixed point at $\tilde{\pi}$, (ii) it is a similarity transformation for each individual contour $c(\tau)$ and (iii) it is a surjection from $\Pi_n$ to $\mathbb{R}^p$.
}

\vspace{0.1in}

As a result of (ii), $h_N^C$ may be viewed as a continuous sequence of similarity mappings from $\mathbb{R}^p$ to $\mathbb{R}^p$ that are indexed by $\tau\in [\tilde{\tau}, +\infty)$. 
The $\tau$-th mapping has expansion factor $\gamma\{N,l(\pi)\}=\gamma(N,\tau)$ and it maps only points on the level-$\tau$ contour $c(\tau)$. Regardless of the amount expanded, $c(\tau)$ and its image are identical in shape. By (\ref{gamma.f}), the expansion factor $\gamma(N,\tau)$ is an increasing function of $\tau$ which approaches infinity when $\tau$ does. Hence, contours farther away from the centre are expanded more and images of the contours fill up the entire $\mathbb{R}^p$. 

If we are to add Condition 5 to Theorem 3, then (iii) can be strengthened to (iii') $h^C_N$ is a bijection from $\Pi_n$ to $\mathbb{R}^p$. See, for example, the proof of Theorem 1 in Tsao and Wu (2013). It is not clear how we may verify condition 5 through $g(X,\theta_x)$ and $g(Y,\theta_y)$. This is why we have not added it to Theorem 3. Nevertheless, we have not encountered any example where Condition 5 is violated. 

\subsection{Extended empirical likelihood on the full parameter space}

By Theorem 3, $h_N^C: \Pi_n\rightarrow \mathbb{R}^p$ is surjective. Thus, for any $\pi \in \mathbb{R}^p$,  $s(\pi)=\{\pi': \pi'\in \Pi_n \hspace{0.05in} \mbox{and} \hspace{0.05in} h_N^C(\pi')=\pi\}$ is non-empty. When $h_N^C$ is not injective, $s(\pi)$ may contain multiple points and $h_N^C$ does not have an inverse. Hence, we define a generalized inverse $h_N^{-C}: \mathbb{R}^p \rightarrow  \Pi_n$ as follows
\begin{equation}
h_N^{-C}(\pi)=\mbox{argmin}_{\pi '\in s(\pi)} \{\|\pi'-\pi\|\}  \hspace{0.2in} \mbox{for $\pi \in \mathbb{R}^p$}. \label{3.29}
\end{equation}
If $s(\pi)$ contains exactly one point $\pi'$, then $h_N^{-C}(\pi)=\pi'$. If $s(\pi)$ has multiple points, then $h_N^{-C}(\pi)$ equals the point $\pi' \in s(\pi)$ that is the closest to $\pi$.

We now define the EEL $l^*(\pi)$ under $h_N^{-C}$ as follows
\begin{equation}
l^*(\pi)=l\{h_N^{-C}(\pi)\} \hspace{0.2in} \mbox{for $\pi \in \mathbb{R}^p$}.  \label{inv.h.oel}
\end{equation}
It is clear that $l^*(\pi)$ is well-defined throughout $\mathbb{R}^p$ since $h_N^{-C}(\pi) \in \Pi_n$ for any $\pi \in \mathbb{R}^p$ and thus the right-hand side of (\ref{inv.h.oel}) is always well-defined. Let $\pi_0'$ be the image of $\pi_0$ under the inversion mapping $h_N^{-C}(\pi)$, that is,
\begin{equation}
h_N^{-C}(\pi_0)=\pi_0'. \label{3.20}
\end{equation}
Then, $l^*(\pi_0)=l\{h_N^{-C}(\pi_0)\}=l(\pi_0')$. Denote by $[\tilde{\pi}, \pi_0]$ the line segment in $\mathbb{R}^p$ that connects the two points $\tilde{\pi}$ and $\pi_0$. Lemma 2 below shows that $\pi_0'$ is on $[\tilde{\pi}, \pi_0]$ and that it is asymptotically very close to $\pi_0$.

\vspace{0.1in}

\noindent {\bf Lemma 2.} \emph{
Under conditions 1, 2 and 3, the point $\pi_0'$ defined by equation (\ref{3.20}) satisfies \\[0.1in]
\hspace*{1.2in} ($i$) $\pi_0' \in [\tilde{\pi}, \pi_0]$ \hspace{0.05in} and \hspace{0.05in} ($ii$) $\pi'_0-\pi_0 = O_p(n^{-3/2})$. %}
}

\vspace{0.1in}

\noindent Theorem 4 below gives the asymptotic distribution of $l^*(\pi_0)$.

\vspace{0.1in}

\noindent {\bf Theorem 4.} \emph{
Under conditions 1, 2 and 3, the two-sample EEL $l^*(\pi)$ defined by (\ref{inv.h.oel})  satisfies
\begin{equation}
l^*(\pi_0) {\longrightarrow} \chi^2_q   \nonumber %\label{3-40}
\end{equation}
in distribution as $n\rightarrow +\infty$.
}

\vspace{0.1in}

The proof of Theorem 4 makes use of the observation that
\begin{equation}
l^*(\pi_0)=l\{h_N^{-C}(\pi_0)\}= l(\pi'_0)=l\{\pi_0+(\pi_0'-\pi_0)\}.  \label{3.35}
\end{equation}
Since  by Lemma 2 $\|\pi_0'-\pi_0\|$ is asymptotically very small, (\ref{3.35}) implies that $l^*(\pi_0)=l(\pi_0)+o_p(1)$. This and the fact that $l(\pi_0)$ has an asymptotic $\chi^2_q$ distribution lead to Theorem 4. The relationship in (\ref{3.35}) is also the key in the derivation of a second-order two-sample EEL in the next section. 

\subsection{Second-order extended empirical likelihood}

We have seen in Theorem 2 that the two-sample OEL admits a Bartlett correction which reduces the coverage error of the empirical likelihood confidence region to $O(n^{-2})$. The following theorem shows that for the just-determined case, the two-sample EEL can also attain the second-order accuracy.

\vspace{0.1in}

\noindent {\bf Theorem 5.} \emph{
 Assume conditions 1, 2, 3 and 4 hold. For the just-determined case where $p=q$, let $l^*_2(\pi)$ be the EEL defined by the composite similarity mapping (\ref{h.function}) with expansion factor $\gamma\{N,l(\pi)\}=\gamma_2\{N,l(\pi)\}$ and
\begin{equation}
\gamma_2\{N,l(\pi)\}=1+\frac{\eta}{2N}\{l(\pi)\}^{\delta(n)}, \label{gamma_2}
\end{equation}
where $\delta(n)=O(n^{-1/2})$ and $\eta$ is the Bartlett correction constant in (\ref{oel.bartfactor}). Then,
\begin{equation}
l^*_2(\pi_0)=l(\pi_0)\{1-\eta N^{-1}+O_p(n^{-3/2})\}. \label{eel2.01}
\end{equation}
and for any fixed $c>0$,
\begin{equation}
P(l_2^*(\pi_0)\leq c)=P(\chi^2_d \leq c)+O(n^{-2}). \label{beel.error}
\end{equation}
}

\vspace{0.1in}

Equation (\ref{beel.error}) follows from (\ref{eel2.01}) and (\ref{oel.bart.00}). It shows that confidence regions based on $l^*_2(\pi)$ have a coverage error of $O(n^{-2})$. Hence, we call $l^*_2(\pi)$ the second-order EEL or EEL$_2$. Correspondingly, we call $l^*(\pi)$ in (\ref{inv.h.oel}), which is defined with the expansion factor $\gamma\{N, l(\pi)\}$ in (\ref{gamma.f}), the first-order EEL or EEL$_1$. The $\delta(n)$ function in $\gamma_2\{N,l(\pi)\}$ is used to control the speed of domain expansion to achieve the second-order accuracy. For convenience, we will use $\delta(n)=n^{-1/2}$ when we compute EEL$_2$ in our numerical examples.

%4444444444444444444444444444444444444444444444444444444444444444444444444444444444
\section{Applications and numerical comparison} 

The need for comparing two populations/models in terms of some numerical aspect of interest arises frequently in applied research. Whenever the numerical aspect of interest can be represented by a parameter defined by estimating equations, the two-sample OEL, BEL and EEL introduced here may be applied to make the comparison. In this section, we consider two such applications. The first is concerned with comparing two populations in terms of the inequality of income distribution. The second is concerned with comparing two linear regression models. Through these two examples, we also compare the numerical accuracy of the three two-sample empirical likelihood methods.

%There are also many other applications of comparing two populations but without strong conditions on their sample data such as simultaneously finding the difference of means and variances.  

%as we know Gini index represents the inequality of income distribution of a country's residents, and the empirical likelihood estimating equation of it is available. Therefore, the application of learning about the difference of the inequality of income distribution from two different countries' residents is applied to our method. There are also many other applications such as simultaneously finding the difference of means and variances of two independent samples or 
%The following two subsections include the results of two applications to illustrate the results of the BEL is always better than OEL, and also show that extended empirical likelihood is always having the best performance among all other empirical likelihood methods.\vspace{0.1in}

\subsection {Application 1: Comparing two Gini indices} 

The Gini index was introduced by Corrado Gini, an Italian statistician of the early 20th century, as a measure of inequality of income or wealth distribution in a country. The value of the Gini index is bounded between 0 and 1, with 0 representing complete equality where all individuals have equal income and 1 representing complete inequality where one individual has all the income and others have none. Gini index has been widely used in social and economic studies of income distributions  [e.g., Gini (1936), Chen (2009) and Domeij et al. (2010), Bee (2012)]. There are also a lot of work on the estimation and inference of the Gini index in both the statistical and econometric literature. 

Qin et al. (2010) and Peng (2011) applied the method of empirical likelihood to make inference about the Gini index. In particular, 
Peng (2011) derived an interesting estimating equation for the Gini index with which the existing theory of Qin and Lawless (1994) was readily applied to make empirical likelihood inference for the index. Peng (2011) also derived empirical likelihoods for the difference between two Gini indexs with paired data and two independent samples. We now apply our two-sample methods to make inference about the difference between two Gini indices using the estimating equation of Peng (2011). For this application, our two-sample OEL coincides with that given by Peng (2011).

Let $X_1,...,X_n$ be i.i.d. observations from an income distribution $F(x)$ supported on $[0, +\infty)$. Define $T_i=\{X_i+X_{[n/2]+i}\}/2$ and $Z_i=\min\{X_i,X_{[n/2]+i}\}$ for $i=1,...,[n/2]$ where $ [ n/2 ]$ is the integer part of $n/2$. Then, Peng (2011) showed that the Gini index, $\theta_0$, of distribution $F(x)$ satisfies,
\begin{equation}
E(T_i-Z_i-T_i\theta_0)=0. \label{ee.for.gini}
\end{equation}
Let $F_A(x)$ be the income distribution of Country A with Gini index $\theta_{x0}$ and $F_B(y)$ be that of Country B with Gini index $\theta_{y0}$.
Suppose we have two random samples of sizes $m$ and $n$, respectively, from $F_A(x)$ and $F_B(y)$. Then, we can compute confidence intervals for the difference $\pi_0=\theta_{y0}-\theta_{x0}$ by using the two-sample OEL, BEL, EEL$_1$ and EEL$_2$. To illustrate their use and to compare the coverage accuracy of confidence intervals based these methods, we consider the following two examples:

\vspace{0.1in}

\hspace*{0.1in} Example 1: $F_A$ is log-normal with $\log(X)\sim N(0,1)$ and $F_B$ is $\chi^2_1$.  
 
\hspace*{0.1in} Example 2: $F_A$ is Pareto$(5)$ and $F_B$ is Exp$(1)$.  \vspace{0.1in}

Before presenting numerical results, note that the two-sample OEL confidence interval for $\pi_0$ is ${\mathcal C}_{1-\alpha}$ given by (\ref{oel.region}). The BEL confidence interval is ${\mathcal C}_{1-\alpha}'$ given in (\ref{bel.region}). 
The EEL$_1$ $l^*(\pi)$ and EEL$_2$ $l_2^*(\pi)$ are both defined through the OEL $l(\pi)$ and the inverse of the composite similarity mapping $h^{-C}_N(\pi)$ in (\ref{inv.h.oel}); the expansion factor in $h^{-C}_N(\pi)$ corresponding to $l^*(\pi)$ is given by (\ref{gamma.f}) and that corresponding to $l_2^*(\pi)$ is given by (\ref{gamma_2}). The EEL$_1$ confidence interval is ${\mathcal C}^*_{1-\alpha}=\{{\pi}: {\pi} \in \mathbb{R}^p \mbox{ and } l^*({\pi}) \leq c\}$ and the EEL$_2$ confidence interval is ${\mathcal C}^{'*}_{1-\alpha}=\{{\pi}: {\pi} \in \mathbb{R}^p \mbox{ and } l^*_2({\pi}) \leq c\}$. The Bartlett correction factor $\eta$ needs to be estimated when computing the BEL and EEL$_2$ confidence intervals, and in both cases we have used the less biased estimator $\tilde{\eta}$ given by Liu and Yu (2010). 

Table 1 contains simulated coverage probabilities of the four confidence intervals for the difference between the two Gini indexes of $F_A$ and $F_B$ in Example 1. Table 2 contains that for Example 2. Each entry in these tables is based on 10,000 pairs of random samples whose sizes are given in the first two columns; it is the proportion of confidence intervals containing the true difference among the 10,000 confidence intervals computed using the 10,000 pairs of samples. We make the following observations based on Tables 1 and 2.

\begin{itemize}

\item[1.] All four confidence intervals give coverage probabilities lower than the nominal level. The OEL interval, in particular, gives the lowest coverage probabilities that may be as much as 10\% lower than the nominal level.  

\item[2.] The BEL, EEL$_1$ and EEL$_2$ intervals are consistently more accurate than the OEL intervals. The two EEL intervals are more accurate than the BEL interval for all combinations of sample sizes and confidence level. Surprisingly, the first-order EEL$_1$ is overall the best, more accurate than the second-order BEL and EEL$_2$ intervals. Hence, we recommend EEL$_1$ for this application.

\item[3.] In column 1 of Table 2, we see that the OEL coverage probability for $(m,n)=(20,40)$ is lower than that for $(20,30)$; in this case the larger sample sizes did not give higher coverage probability. This surprising phenomenon occurs sometimes for other other two-sample methods as well. See also Table 2 in Liu and Yu (2010) for similar results. Noting that $m-n$ is smaller in $(20, 30)$, it seems that for two-sample inference a large difference in sample size can negatively affect the accuracy of the EL based confidence intervals.

\end{itemize}  
\subsection {Application 2: Comparing two linear regression models} Consider two simple linear regression models having the same predictor variable but possibly different slopes and intercepts. To compare the parameters of the two models with two independent random samples (one from each model), a commonly used method is to introduce a dummy/indicator variable and the comparison is then done through a multiple linear regression model with two covariates; the predictor variable and the dummy variable. This method, however, requires the assumption that error distributions of the two models are the same. Without making this assumption, we now use two-sample empirical likelihood methods to compare the model parameters. Specifically, we compare models
\[(a) \hspace{0.1in} y={x}^T\beta_a +\varepsilon_a \hspace{0.2in} \mbox{and} \hspace{0.2in} (b) \hspace{0.1in} y={x}^T\beta_b +\varepsilon_b, \]
where $\beta_a=(\beta_{a0},\beta_{a1})^T$, $\beta_b=(\beta_{b0},\beta_{b1})^T$, $\varepsilon_a$ and $\varepsilon_b$ are random errors with possibly different distributions,  but ${x}=(1,x_1)^T$ is the same in both models. The parameter vector of interest is the difference $\pi=\beta_a-\beta_b$.

For our simulation study, $x_1$ is assumed to be a uniform random variable on $[0,30]$. We consider the following two examples:

\vspace{0.1in}

\hspace*{0.1in} Example 3: Model (a) with $\varepsilon_a \sim N(0,1)$ and ${\beta_a}=(2,1)^T$ and Model (b) with $\varepsilon_b \sim N(0,1)$ and ${\beta_b}=(2,2)^T$. 
 
\hspace*{0.1in} Example 4: Model (a) with $\varepsilon_a \sim Exp(1)-1$ and ${\beta_a}=(2,1)^T$ and Model (b) with $\varepsilon_b \sim N(0,1)$ and ${\beta_b}=(2,2)^T$. 

The simulated coverage probabilities for $\pi$ given by the four empirical likelihood methods are shown in Table 3. Although Examples 3 and 4 are multi-dimensional examples ($p=2$), the three observations made above also apply to Table 3. In particular, overall EEL$_1$ has better accuracy than OEL, BEL, and EEL$_2$. We recommend EEL$_1$ due to its simplicity and excellent accuracy.

%\subsection{Applications results}

%We draw 10,000 random samples for each case of above two applications to illustrate the coverage accuracy of the two-sample empirical likelihood for estimating equations. The applications results show that the BEL is always more accurate than the OEL which proved the theoretical finding. 

\begin{table}

\scriptsize{
\title{Coverage probabilities ($\%$) of confidence regions based on OEL, EEL$_1$, BEL and EEL$_2$ for Example 1}{%
\begin{tabular}{lllcccccccccccc} \\
& &  & \multicolumn{4}{c}{90\% level} & \multicolumn{4}{c}{95\% level} & \multicolumn{4}{c}{99\% level}\\
&$m$&    $n$&OEL&EEL$_1$&BEL &EEL$_2$ &OEL&EEL$_1$&BEL &EEL$_2$ &OEL&EEL$_1$&BEL &EEL$_2$\\ \hline
&20   &20	&80.0	&81.9	&81.4	&82.4	&86.5	&88.7	&87.7	&88.4	&94.0	&95.8	&94.5	&95.2	\\
&     &30	&81.1	&83.3	&82.4	&83.9	&87.6	&89.8	&88.8	&89.8	&95.2	&96.8	&95.7	&96.3	\\
&     &40	&82.0	&84.1	&83.2	&84.9	&88.5	&90.9	&89.6	&91.0	&95.3	&96.9	&95.7	&96.3	\\
&     &60	&82.1	&84.2	&83.1	&85.7	&88.1	&90.7	&89.2	&91.3	&95.9	&97.2	&96.4	&97.0	\\
&30	  &20	&79.7	&81.7	&80.9	&82.2	&86.6	&88.4	&87.5	&88.5	&94.0	&95.8	&94.6	&95.2	\\
&	    &30	&82.6	&84.0	&83.7	&84.5	&89.1	&90.3	&89.9	&90.3	&95.7	&96.6	&96.0	&96.3	\\
&	    &40	&83.2	&84.5	&84.2	&85.1	&89.5	&90.9	&90.3	&91.0	&96.0	&97.1	&96.4	&96.9	\\
&	    &60	&84.2	&85.5	&84.9	&86.3	&90.3	&91.7	&91.1	&92.0	&97.0	&97.8	&97.3	&97.7	\\
&40 	&20	&80.4	&82.2	&81.5	&82.9	&87.0	&88.5	&87.8	&88.6	&94.2	&95.7	&94.7	&95.2	\\
&	    &30	&82.9	&84.2	&83.8	&84.8	&89.6	&90.9	&90.4	&90.9	&96.2	&97.2	&96.6	&96.9	\\
&	    &40	&84.4	&85.5	&85.2	&86.0	&90.6	&91.6	&91.3	&91.7	&96.8	&97.6	&97.1	&97.4	\\
&	    &60	&85.6	&86.6	&86.5	&87.2	&91.4	&92.5	&92.1	&92.8	&97.2	&97.8	&97.5	&97.7	\\
&60	  &20	&79.7	&81.5	&80.8	&82.7	&86.7	&88.6	&87.5	&88.8	&94.6	&95.8	&95.0	&95.6	\\
&    	&30	&83.5	&84.7	&84.2	&85.3	&89.6	&90.6	&90.2	&91.0	&96.0	&97.0	&96.3	&96.8	\\
&    	&40	&85.1	&86.0	&85.9	&86.7	&91.3	&92.1	&91.9	&92.3	&97.2	&97.8	&97.4	&97.6	\\
&    	&60	&85.8	&86.4	&86.5	&86.9	&91.6	&92.4	&92.2	&92.5	&97.4	&97.8	&97.5	&97.7	\\\hline
\end{tabular}}
}
\label{table1}
%\begin{tabnote}
\noindent Each entry in the table is a simulated coverage probability for ${\pi}$ based on 10,000 random samples of size $m$ and $n$ indicated in column 1 and 2 from the distribution log-normal (i.e. log $N(0,1)$) and $\chi^2_1$, respectively.
%\end{tabnote}
\end{table}

\begin{table}
\scriptsize{
%\tbl{Coverage probabilities of confidence regions based on the original empirical likelihood (OEL), the first-order extended empirical likelihood (EEL$_1$), the Bartlett corrected empirical likelihood (BEL) and Bartlett corrected extended empirical likelihood (EEL$_2$) (pareto5-exp1)}{%
\title{Coverage probabilities ($\%$) of confidence regions based on OEL, EEL$_1$, BEL and EEL$_2$ for Example 2}{%
\begin{tabular}{lllcccccccccccc} \\
& &  & \multicolumn{4}{c}{90\% level} & \multicolumn{4}{c}{95\% level} & \multicolumn{4}{c}{99\% level}\\
&$m$&    $n$&OEL&EEL$_1$&BEL &EEL$_2$ &OEL&EEL$_1$&BEL &EEL$_2$ &OEL&EEL$_1$&BEL &EEL$_2$\\ \hline
&20 &20&80.8	&83.4	&81.8	&83.4	&86.8	&89.7	&87.6	&89.3	&93.6	&96.5	&94.1	&95.8	\\
&	  &30	&82.0	&84.5	&83.1	&84.8	&88.1	&90.1	&88.7	&89.9	&94.3	&96.2	&94.7	&95.7	\\
&	  &40	&80.1	&84.4	&81.0	&84.7	&85.6	&90.4	&86.4	&90.2	&91.7	&96.2	&92.2	&95.9	\\
&	  &60	&82.1	&83.8	&83.0	&84.7	&88.0	&89.8	&88.8	&90.1	&95.0	&96.4	&95.3	&96.1	\\
&30	&20	&84.1	&86.0	&85.1	&86.2	&90.1	&92.0	&90.7	&91.7	&95.8	&97.7	&96.2	&97.2	\\
&	  &30	&84.2	&86.1	&85.0	&86.3	&89.9	&91.9	&90.4	&91.7	&95.8	&97.5	&96.1	&97.2	\\
&	  &40	&84.2	&86.2	&85.0	&86.4	&89.9	&92.2	&90.6	&92.1	&95.2	&97.5	&95.5	&97.1	\\
&	  &60	&84.1	&85.4	&84.9	&86.0	&90.2	&91.6	&90.8	&91.7	&96.2	&97.2	&96.5	&97.0	\\
&40	&20	&85.3	&87.3	&86.2	&87.5	&91.1	&93.0	&91.6	&92.9	&96.6	&98.3	&96.8	&97.9	\\
&   & 30&85.1	&87.2	&85.8	&87.5	&90.8	&92.9	&91.2	&92.8	&96.4	&98.3	&96.7	&98.0	\\
&   & 40&86.0	&87.1	&86.7	&87.3	&91.3	&92.5	&91.8	&92.5	&97.0	&98.0	&97.1	&97.7	\\
&   & 60&85.4	&86.4	&86.1	&86.7	&91.2	&92.1	&91.7	&92.2	&96.8	&97.6	&97.1	&97.4	\\
&60	&20	&86.2	&88.0	&86.8	&88.7	&91.8	&93.9	&92.3	&93.8	&97.6	&99.0	&97.7	&98.7	\\
&   & 30&86.4	&88.6	&87.0	&89.0	&91.8	&94.1	&92.2	&94.1	&96.9	&98.6	&97.1	&98.4	\\
&   & 40&87.2	&88.4	&87.8	&88.7	&92.8	&93.9	&93.3	&93.9	&98.0	&98.6	&98.1	&98.4	\\
&   &	60&86.9	&88.0	&87.5	&88.2	&92.5	&93.3	&92.7	&93.3	&97.4	&98.4	&97.6	&98.2	\\\hline
\end{tabular}}
}
\label{table2}
%\begin{tabnote}
\noindent Each entry in the table is a simulated coverage probability for ${\pi}$ based on 10,000 random samples of size $m$ and $n$ indicated in column 1 and 2 from the distribution $Pareto(5)$ and $Exp(1)$, respectively.
%Each entry in the table is a simulated coverage probability for ${\beta}$ based on 10,000 random samples of size $n$ indicated in column 2 from the linear model indicated in column 1.
%\end{tabnote}
\end{table}

\begin{table}
\scriptsize{
%\tbl{Coverage probabilities of confidence regions based on the original empirical likelihood (OEL), the first-order extended empirical likelihood (EEL$_1$), the Bartlett corrected empirical likelihood (BEL) and Bartlett corrected extended empirical likelihood (EEL$_2$) (regression)}{%
\title{Coverage probabilities ($\%$) of confidence regions based on OEL, EEL$_1$, BEL and EEL$_2$ for Example 3 (Ex-3) and Example 4 (Ex-4)}{%
\begin{tabular}{llcccccccccccc} \\
  & & \multicolumn{4}{c}{90\% level} & \multicolumn{4}{c}{95\% level} & \multicolumn{4}{c}{99\% level}\\
&($m$,$n$)&OEL&EEL$_1$&BEL &EEL$_2$ &OEL&EEL$_1$&BEL &EEL$_2$ &OEL&EEL$_1$&BEL &EEL$_2$\\ \hline
{Ex-3} 	&(20, 20)	&80.9	&85.8	&83.0	&85.2	&86.9	&91.6	&88.6	&90.5	&93.2	&96.9	&94.1	&95.8	\\
     	&(20, 40)	&81.2	&85.2	&83.0	&84.7	&86.8	&90.7	&88.1	&89.8	&93.2	&96.6	&93.8	&95.5	\\
     	&(30, 30)	&84.2	&87.5	&85.8	&87.3	&90.2	&93.2	&91.5	&92.9	&95.4	&97.7	&95.7	&97.1	\\
     	&(40, 30)	&83.1	&87.0	&84.2	&86.5	&88.4	&92.6	&89.5	&91.7	&93.9	&97.2	&94.3	&96.4	\\
     	&(40, 40)	&87.2	&88.6	&88.5	&88.6	&92.9	&94.2	&93.8	&94.0	&98.1	&98.9	&98.6	&98.7	\\
     	&(50, 30)	&80.7	&86.3	&81.9	&85.5	&86.1	&91.5	&87.0	&90.7	&91.3	&96.0	&91.6	&95.1	\\
     	&(50, 50)	&87.7	&88.6	&88.5	&88.7	&93.6	&94.5	&94.2	&94.3	&98.5	&99.0	&98.8	&98.8	\\
{Ex-4}&(20, 20)	&78.0	&83.3	&80.6	&83.0	&84.4	&89.5	&86.4	&88.7	&91.5	&95.6	&92.6	&94.6	\\
     	&(20, 40)	&79.4	&84.0	&81.4	&83.5	&85.3	&89.4	&86.5	&88.6	&91.6	&95.2	&92.2	&94.3	\\
     	&(30, 30)	&82.7	&86.4	&84.6	&86.3	&88.7	&92.1	&90.0	&91.7	&94.4	&97.2	&95.0	&96.5	\\
     	&(40, 30)	&80.3	&85.1	&81.8	&84.8	&86.4	&91.2	&87.6	&90.5	&91.8	&95.9	&92.3	&95.1	\\
      &(40, 40)	&86.1	&87.7	&87.6	&88.0	&92.4	&93.6	&93.4	&93.6	&97/7	&98.5	&98.2	&98.3	\\
     	&(50, 30)	&78.7	&85.0	&80.0	&84.1	&84.2	&90.3	&85.3	&89.3	&89.7	&95.0	&90.1	&94.1	\\
     	&(50, 50)	&86.7	&87.8	&87.8	&88.2	&92.2	&93.4	&93.2	&93.4	&97.7	&98.5	&98.2	&98.3	\\  \hline
\end{tabular}}
}
\label{table3}
%\begin{tabnote}
\noindent Each entry in the table is a simulated coverage probability for ${\pi}$ based on 10,000 pairs of random samples with sizes $(m,n)$ indicated in column 2 from the linear models indicated in column 1.
%\end{tabnote}
\end{table}

%5555555555555555555555555555555555555555555555555555555555555555555
\section{Appendix}

We now present proofs of theorems and lemmas in the order as they appeared in the paper. For brevity, for results that are minor variations of existing results in the literature, we give only references to the existing results instead of detailed proofs which may be found in the references.

\vspace{0.1in}

\noindent Theorem 1 is the standard first-order result for an OEL. It is implied by Theorem 2 which gives the second-order result. Hence, its proof is omitted. We now prove Lemma 1 by following that for equation (3) in Liu et al. (2008).

\vspace{0.1in}

\noindent {\bf Proof of Lemma 1}\\[0.1in]
First note that, under conditions of Lemma 1, $\lambda^*=O_p(n^{-1/2})$; see the proof of Theorem 1 in Owen (1990). For clarity of presentation, we break the proof of Lemma 1 into the following three steps.
%For simplicity, from now on, we use the notations $\theta_x^{0}$ and $\theta_x^{1}$ to represent $\theta^0_{y}-\pi^0$ and $\theta^1_{y}-\pi^1$, respectively.\\

\vspace{0.1in}

\noindent Step 1: Let $C_{11}=\frac{1}{n}\sum g(y_j,\theta_y^{0})$ and $C_{12}=\frac{1}{m}\sum g(x_i,\theta_x^{0})$. Taylor expansion of the first equation of (\ref{oel.xy}) gives %(\ref{oel.log}), we have
\begin{equation}
	\frac{1}{n}\sum g(y_j,\theta^*_{y}) -\left\{\frac{f_n}{n} \sum g(y_j,\theta^*_{y})g^T(y_j,\theta^*_{y}) - V_1 \right\} \lambda^*-V_1\lambda^* +O_p(n^{-1})=0. \label{oel.ta1}
\end{equation} 
It follows that
\begin{equation}
	\lambda^*=V_1^{-1}\left\{\frac{1}{n}\sum g(y_j,\theta_y^{0})\right\} +O_p(n^{-1})=V_1^{-1}C_{11} +O_p(n^{-1}). \label{oel.ta2}
\end{equation} 
Similarly, expansion of the second equation of (\ref{oel.xy}) gives,
\begin{equation}
	\lambda^*=-V_2^{-1}\left\{\frac{1}{m}\sum g(x_i,\theta_x^{0})\right\} +O_p(n^{-1})=-V_2^{-1}C_{12}+O_p(n^{-1}). \label{oel.ta5}
\end{equation} 
Based on (\ref{oel.ta2}) and (\ref{oel.ta5}),
\begin{equation}\nonumber 
	V_1\lambda^*=V_1^{-1}C_{11} +O_p(n^{-1}). 
\end{equation} 
\begin{equation}\nonumber 
	-V_2^{-1}\lambda^*=C_{12}+O_p(n^{-1}). 
\end{equation} 
thus, we have
\[\zeta^0=V^{-1}(C_{11}-C_{12})= V^{-1}D_1. \]
Step 2.  Further expanding the left-hand side of (\ref{oel.ta1}), we have
\begin{eqnarray}\label{oel.ta3}
	&&\frac{1}{n}\sum g(y_j,\theta^*_{y}) -\left\{\frac{f_n}{n} \sum g(y_j,\theta_y^{0})g^T(y_j,\theta_y^{0}) - V_1 \right\}\zeta^0 \\ 
	&&-V_1\lambda^* +\frac{f_n^2}{n} \sum \left\{(\zeta^0)^Tg(y_j,\theta_y^{0})\right\}^2g(y_j,\theta_y^{0}) +O_p(n^{-3/2})=0.  \nonumber
\end{eqnarray} 
Let 
\begin{eqnarray}
	&&C_{21}=-\left\{\frac{f_n}{n} \sum g(y_j,\theta_y^{0})g^T(y_j,\theta_y^{0}) - V_1 \right\}\zeta^0 ,\\ \nonumber 
	&&C_{31}= \frac{f_n^2}{n} \sum \left\{(\zeta^0)^Tg(y_j,\theta_y^{0})\right\}^2g(y_j,\theta_y^{0}),\\ \nonumber  
	&&C_{22}=\left\{ \frac{f_m}{m} \sum g(x_i,\theta_x^{0})g^T(x_i,\theta_x^{0}) - V_2 \right\}\zeta^0 , \\ \nonumber 
	&&C_{32}= -\frac{f_m^2}{m} \sum \left\{(\zeta^0)^Tg(x_i,\theta_x^{0})\right\}^2g(x_i,\theta_x^{0}) \nonumber 
\end{eqnarray}
It follows that
\begin{eqnarray} \nonumber
	&&\lambda^*=V_1^{-1}(C_{11}+C_{21}+C_{31}) +O_p(n^{-3/2}), \\ \nonumber
	&&\lambda^*=-V_2^{-1}(C_{12}+C_{22}+C_{32}) +O_p(n^{-3/2}), \nonumber
\end{eqnarray} 
and we obtain  
\begin{equation}
	\zeta^1=V^{-1}\left\{(C_{11}-C_{12})+(C_{21}-C_{22})+(C_{31}-C_{32})\right\}= V^{-1}(D_1+D_2+D_3). \nonumber
\end{equation} 
Step 3. Additional expanding (\ref{oel.ta3}) gives
\begin{eqnarray} \nonumber
	&&\frac{1}{n}\sum g(y_j,\theta^*_{y}) -\left\{\frac{f_n}{n} \sum g(y_j,\theta_y^{1})g^T(y_j,\theta_y^{1}) - V_1 \right\}\zeta^1  -V_1\lambda^* +\\ \nonumber
	&&\frac{f_n^2}{n} \sum \left\{(\zeta^1)^Tg(y_j,\theta_y^{1})\right\}^2g(y_j,\theta_y^{1})+ \frac{f_n^3}{n} \sum \left\{(\zeta^0)^Tg(y_j,\theta_y^{0})\right\}^3g(y_j,\theta_y^{0})+O_p(n^{-2})=0, \label{oel.ta4}
\end{eqnarray} 
Let 
\begin{eqnarray} \nonumber 
	&&C^*_{21}=-\left\{ \frac{f_n}{n} \sum g(y_j,\theta_y^{1})g^T(y_j,\theta_y^{1}) - V_1 \right\}\zeta^1 ,\\ \nonumber 
	&&C^*_{31}=\frac{f_n^2}{n} \sum \left\{(\zeta^1)^Tg(y_j,\theta_y^{1})\right\}^2g(y_j,\theta_y^{1}),\\ \nonumber  
	&&C_{41}=-\frac{f_n^3}{n} \sum \left\{(\zeta^0)^Tg(y_j,\theta_y^{0})\right\}^3g(y_j,\theta_y^{0}),\\ \nonumber  
	&&C^*_{22}=\left\{ \frac{f_m}{m} \sum g(x_i,\theta_x^{1})g^T(x_i,\theta_x^{1}) - V_2 \right\}\zeta^1 , \\ \nonumber 
	&&C^*_{32}= -\frac{f_m^2}{m} \sum \left\{(\zeta^1)^Tg(x_i,\theta_x^{1})\right\}^2g(x_i,\theta_x^{1}), \\\nonumber 
	&&C_{42}= \frac{f_m^3}{m} \sum \left\{(\zeta^0)^Tg(x_i,\theta_x^{0})\right\}^3g(x_i,\theta_x^{0}). \nonumber 
\end{eqnarray}
which gives the following higher order expansions of $\lambda^*$,
\begin{eqnarray}
	&&\lambda^*=V_1^{-1}(C_{11}+C^*_{21}+C^*_{31}+C_{41}) +O_p(n^{-2}), \\ \nonumber
	&&\lambda^*=-V_2^{-1}(C_{12}+C^*_{22}+C^*_{32}+C_{42}) +O_p(n^{-2}).  \nonumber
\end{eqnarray} 
Thus
\begin{eqnarray} 
	\zeta^2&=&V^{-1}\left\{(C_{11}-C_{12})+(C^*_{21}-C^*_{22})+(C^*_{31}-C^*_{32})+(C_{41}-C_{42})\right\} \\ \nonumber
	&=& V^{-1}(D_1+D^*_2+D^*_3+D_4).  \nonumber
\end{eqnarray} 
Then, the Taylor expansion for $l(\pi)/N$ can be expressed as
\begin{eqnarray}
	\frac{l(\pi)}{N}=&&\frac{2}{n}\sum (\zeta^2)^T g(y_j,\theta_y^0) - \frac{2}{m}\sum (\zeta^2)^T g(x_i,\theta_x^0) \\ \nonumber
	&&-\frac{f_n}{n}\sum (\zeta^2)^T g(y_j,\theta_y^1)g(y_j,\theta_y^1)^T\zeta_2 - \frac{f_m}{m}\sum (\zeta^2)^Tg(x_i,\theta_x^1)g(x_i,\theta_x^1)^T\zeta_2 \\ \nonumber
	&&+\frac{2f_n^2}{3n}\sum \left\{(\zeta^1)^T g(y_j,\theta_y^1)\right\}^3 - \frac{2f_m^2}{3m}\sum \left\{(\zeta^1)^Tg(x_i,\theta_x^1)\right\}^3 \\ \nonumber
	&&+\frac{f_n^3}{2n}\sum \left\{(\zeta^0)^T g(y_j,\theta_y^0)\right\}^4 - \frac{f_m^3}{2m}\sum \left\{(\zeta^0)^Tg(x_i,\theta_x^0)\right\}^4+O_p(n^{-5/2}) \\ \nonumber
	=&&2I_1-I_2+\frac{2}{3}I_3-\frac{1}{2}I_4+O_p(n^{-5/2})
\end{eqnarray}
where 
\begin{eqnarray}
	I_1&=&\zeta^2(C_{11}-C_{12})=D_1^TV^{-1}(D_1+D^*_2+D^*_3+D_4), \\ \nonumber
	I_2&=&\zeta^2(C_{11}-C_{12}+C^*_{31}-C^*_{32}+C_{41}-C_{42})=(D_1+D^*_3+D_4)V^{-1}(D_1+D^*_2+D^*_3+D_4), \\ \nonumber
	I_3&=&\zeta^1(C^*_{31}-C^*_{32})=(D^*_3)V^{-1}(D_1+D_2+D_3), \\ \nonumber
	I_4&=&-\zeta^0(C_{41}-C_{42})=-D_1V^{-1}(D_4). \nonumber
\end{eqnarray}
Noting that 
\begin{eqnarray}
	\nonumber &&D_1=O_p(n^{-1/2}),  \hspace{0.1in} D_2=O_p(n^{-1}), \hspace{0.1in} D^*_2=O_p(n^{-1}),\hspace{0.1in} D_2-D^*_2=O_p(n^{-3/2}), \\ \nonumber 
	&&D_3=O_p(n^{-1}),  \hspace{0.1in} D^*_3=O_p(n^{-1}), \hspace{0.1in} D_3-D^*_3=O_p(n^{-3/2}),\hspace{0.1in} D_4=O_p(n^{-3/2}),  \nonumber 
\end{eqnarray}
thus, we have 
\begin{eqnarray}
	\frac{l(\pi)}{N} \nonumber 
	=&&2D_1^TV^{-1}(D_1+D^*_2+D^*_3+D_4)- (D_1+D^*_3+D_4)V^{-1}(D_1+D^*_2+D^*_3+D_4) \\ \nonumber
	&&+\frac{2}{3}(D^*_3)V^{-1}(D_1+D_2+D_3)+\frac{1}{2}D_1V^{-1}D_4+ O_p(n^{-5/2})\\ \nonumber
	=&&  D_1^TV^{-1}(D_1+D^*_2+D^*_3+D_4)- (D^*_3+D_4)V^{-1}(D_1+D^*_2+D^*_3+D_4)     \\ \nonumber
	&&+\frac{2}{3}(D^*_3)V^{-1}(D_1+D_2+D_3)+\frac{1}{2}D_1V^{-1}D_4+ O_p(n^{-5/2}) \\ \nonumber                           
	=&&\left\{D_1^TV^{-1}D_1+D_1^TV^{-1}D^*_2+D_1^TV^{-1}D^*_3+D_1^TV^{-1}D_4\right\} \\ \nonumber
	&&-\left\{D^{*T}_3V^{-1}D_1+D^{*T}_3V^{-1}D^*_2+D^{*T}_3V^{-1}D^*_3+D_4^TV^{-1}D_1+O_p(n^{-5/2})\right\} \\ \nonumber
	&&+\frac{2}{3}\left\{D^{*T}_3V^{-1}D_1+D^{*T}_3V^{-1}D_2+D^{*T}_3V^{-1}D_3\right\} +\frac{1}{2}D_1^TV^{-1}D_4+O_p(n^{-5/2}) \\ \nonumber
	=&& D_1^TV^{-1}D_1+D_1^TV^{-1}D^*_2+\frac{2}{3}D^{T}_1V^{-1}D_3^*-\frac{1}{3}D^{T}_2V^{-1}D_3\\ \nonumber
	&&-\frac{1}{3}D^{T}_3V^{-1}D_3+\frac{1}{2}D_1^TV^{-1}D_4+O_p(n^{-5/2}). \nonumber
\end{eqnarray}
Hence,
\begin{eqnarray}
	\frac{l(\pi)}{N} \nonumber &=&S^{\tau}S^{\tau}-S^{\tau\upsilon}S^{\tau}S^{\upsilon}+\frac{2}{3}s^{\tau\alpha\beta}S^{\tau}S^{\alpha}S^{\beta}+S^{\tau\upsilon}S^{\upsilon\omega}S^{\tau}S^{\omega}+\frac{2}{3}S^{\tau\alpha\beta}S^{\tau}S^{\alpha}S^{\beta},\\  \nonumber       &&-2s^{\tau\upsilon\omega}S^{\tau\alpha}S^{\upsilon}S^{\alpha}S^{\upsilon}+s^{\tau\upsilon\omega}s^{\tau\alpha\beta}S^{\upsilon}S^{\omega}S^{\alpha}S^{\beta}-\frac{1}{2}s^{\tau\upsilon\omega\alpha}S^{\tau}S^{\upsilon}S^{\omega}S^{\alpha} \\ \nonumber
	&& +(S^{\tau\upsilon}-S_1^{\tau\upsilon})S^{\tau}S^{\upsilon}+\frac{2}{3}(S_1^{\tau\alpha\beta}-S^{\tau\alpha\beta})S^{\tau}S^{\alpha}S^{\beta} +O_p(n^{-5/2}). \nonumber
\end{eqnarray}
where 
\begin{eqnarray}
	\nonumber&& D_1^TV^{-1}D_1 =S^{\tau}S^{\tau}, \hspace{0.1in} \\\nonumber
	&& D_1^TV^{-1}D^*_2=-S_1^{\tau\upsilon}S^{\tau}S^{\upsilon}+ S^{\tau\upsilon}S^{\upsilon\omega}S^{\tau}S^{\omega}- s^{\omega\alpha\beta}S^{\tau\upsilon}S^{\upsilon\omega}S^{\alpha}S^{\beta}+O_p(n^{-5/2})\\ \nonumber 
	&& D^{T}_1V^{-1}D_3^*=(S_1^{\tau\alpha\beta}-S^{\tau\alpha\beta})S^{\tau}S^{\alpha}S^{\beta} ,\hspace{0.1in} D^{T}_2V^{-1}D_3=-s^{\tau\upsilon\omega}S^{\tau\alpha}S^{\upsilon}S^{\alpha}S^{\upsilon} \\ \nonumber
	&& D^{T}_3V^{-1}D_3=s^{\tau\upsilon\omega}s^{\tau\alpha\beta}S^{\upsilon}S^{\omega}S^{\alpha}S^{\beta}, \hspace{0.1in} D_1^TV^{-1}D_4=s^{\tau\upsilon\omega\alpha}S^{\tau}S^{\upsilon}S^{\omega}S^{\alpha} \nonumber
\end{eqnarray}
which proves Lemma 1. 
%\end{proof}
%\begin{proof}{Proof of Theorem 3}{}

\vspace{0.1in}

\noindent With Lemma 1, the proof of Theorem 2 follows from that for the second order result in DiCiccio et al. (1991). See DiCiccio et al. (1988) for details.

\vspace{0.1in}

%\end{proof}
%%%%%%%%%%%%%%%%%%%%%%%%%%%%%%%%%%%%%%%%%%%%%%%%%%%%%%%%%%%%%%%%%%%%%%%%%%%%%%%%%%%%%%%%%%%%%%%%%%%%%%%%%%%
\noindent {\bf Proof of Lemma 2}\\[0.1in]
Differentiating $l(\pi)$ in (2) and evaluating the derivative at $\pi_0$, we find $J(\pi_0) =\frac{\partial l(\pi)}{\partial \pi}|_{\pi=\pi_0}$ as follows
%\begin{footnotesize}
%\begin{eqnarray}
%\nonumber J(\pi_0)  &=&\frac{\partial l(\pi)}{\partial \pi}|_{\pi=\pi_0} \\ 
% &=&2\lambda^T(\pi_0)\left\{ \sum 
%              \frac{f_n g'(y_j,\theta_{y_0})\frac{\partial \theta_y}{\partial \pi}|_{\pi=\pi_0} }{1+f_n \lambda^T(\pi_0) g(y_j,\theta_{y_0})}
%              -\sum \frac{f_m g'(x_i,\theta_{x_0})\frac{\partial \theta_x}{\partial \pi}|_{\pi=\pi_0}}{1-f_m \lambda^T(\pi_0)g(x_i,\theta_{x_0})} \right\},           %\label{lemma1.01}
%\end{eqnarray}
%\end{footnotesize}

\begin{equation}
J(\pi_0)  =\lambda^T(\pi_0)\left\{ \sum 
              \frac{f_n g'(y_j,\theta_{y_0})\frac{\partial \theta_y}{\partial \pi}|_{\pi=\pi_0} }{1+f_n \lambda^T(\pi_0) g(y_j,\theta_{y_0})}
              -\sum \frac{f_m g'(x_i,\theta_{x_0})\frac{\partial \theta_x}{\partial \pi}|_{\pi=\pi_0}}{1-f_m \lambda^T(\pi_0)g(x_i,\theta_{x_0})} \right\},           \label{lemma1.01}
\end{equation}

\noindent where $g'(y_j,\theta_{y_0})\frac{\partial \theta_y}{\partial \pi}|_{\pi=\pi_0} =\partial g(y_j,\theta_y)/\partial \pi|_{\pi=\pi_0}$ and $g'(x_i,\theta_{x_0})\frac{\partial \theta_x}{\partial \pi}|_{\pi=\pi_0}=\partial g(x_i,\theta_x)/\partial \pi|_{\pi=\pi_0}$.  Under the conditions of the lemma, we can show that 
$\lambda(\pi_0)=O_p(n^{-1/2})$ and $J(\pi_0)=O_p(n^{1/2})$.
% $\lambda(\pi_0)=O_p(n^{-1/2})$, $J(\pi_0)=O_p(n^{1/2})$ and $[\partial^2 l(\pi)/\partial \pi^2]_{\pi=\pi_0}=O_p(n)$.
Also, applying Taylor expansion to $l(\pi)$ in a small neighbourhood of $\pi_0$, $\{\pi: \|\pi-\pi_0\|\leq \kappa n^{-1/2}\}$, where $\kappa$ is some positive constant, we obtain
\begin{equation}
 l(\pi)=l\{\pi_0+(\pi-\pi_0)\}=l(\pi_0)+J(\pi_0)(\pi-\pi_0) + O_p(1). \label{lemma1.02}
 \end{equation}
By Theorem 1, $l(\pi_0)=O_p(1)$. This and (\ref{lemma1.02}) imply that for a $\pi$ in that small neighbourhood,
\begin{equation}
l(\pi)=O_p(1). \label{lemma1.03}
\end{equation}
To show part (i), since $h_N^C(\pi_0')=\pi_0$, we have
\begin{equation}  \pi_0-\tilde{\pi}=\gamma\{n,l(\pi_0')\}(\pi_0'-\tilde{\pi}). \label{lemma1.04} \end{equation}
Noting that $\gamma\{N,l(\pi)\}\geq 1$, (\ref{lemma1.04}) implies that $\pi'_0$ is on the ray originating from $ \tilde{\pi}$ through $\pi_0$ and $$\|\pi_0-\tilde{\pi}\|\geq\|\pi_0'-\tilde{\pi}\|.$$ 
Hence, $\pi_0'\in [\tilde{\pi}, \pi_0]$ and part (i) of the lemma 2 is proven.

To show part (ii), since $\tilde{\pi}$ is $\surd{n}$-consistent and $\pi_0'\in [\tilde{\pi}, \pi_0]$, we have $\pi_0'-\pi_0 =O_p(n^{-1/2})$. It follows from (\ref{lemma1.03}) that $l(\pi_0')=O_p(1)$. This implies
\begin{equation}
\gamma\{N,l(\pi_0')\}=1+\frac{l(\pi_0')}{2N}=1+O_p(n^{-1}). \label{lemma1.05}
\end{equation}
Adding and subtracting a $\pi_0$ on the right-hand side of (\ref{lemma1.04}) gives
\begin{equation} \pi_0-\tilde{\pi}=\gamma\{N,l(\pi_0')\}(\pi_0'-\pi_0+\pi_0-\tilde{\pi}). \nonumber \end{equation}
This implies that
\begin{equation}
 \left[1-\gamma\{N,l(\pi_0')\}\right](\pi_0-\tilde{\pi})=\gamma\{N,l(\pi_0')\}(\pi_0'-\pi_0).  \label{lemma1.06}
\end{equation}
It follows from (\ref{lemma1.05}), (\ref{lemma1.06}) and $\tilde{\pi}-\pi_0=O_p(n^{-1/2})$ that 
\begin{equation}
\pi_0'-\pi_0=O_p(n^{-3/2}). \nonumber %\label{b-15}
\end{equation}
This proves part (ii) of the lemma 2.
$\Box$ \vspace{0.1in}

Proof for Theorem 3 follows easily from that for Theorem 1 in Tsao and Wu (2014).

%******************* Theorem 4

\noindent {\bf Proof of Theorem 4}\\[0.1in]
By (ii) of Lemma 2, $\pi_0'-\pi_0=O_p(n^{-3/2})$. Taylor expansion of $l^*(\pi_0)$ gives
\begin{equation}
 l^*(\pi_0)= l(\pi_0')=l\{\pi_0+(\pi_0'-\pi_0)\}=l(\pi_0)+J(\pi_0)(\pi_0'-\pi_0) + o_p(n^{-3/2}). \label{t2-01}
 \end{equation}
Since $J(\pi_0)=O_p(n^{1/2})$, (\ref{t2-01}) implies that 
$l^*(\pi_0)=l(\pi_0)+ O_p(n^{-1})$. Thus, the extended empirical log-likelihood ratio $l^*(\pi_0)$ has the same limiting $\chi^2_q$ distribution as the original empirical log-likelihood ratio $l(\pi_0)$. 
$\Box$ \vspace{0.1in}

We need the following lemma for the proof of Theorem 5. \vspace{0.1in}

\noindent \textbf{Lemma 3.} {\em Assume conditions 1, 2, 3 and 4 hold. Under the composite similarity mapping (9) with expansion factor $\gamma\{N,l(\pi)\}=\gamma_2\{N, l(\pi)\}$ in (17), we have}
\begin{equation}
\pi_0'-\pi_0=\frac{b}{2n}(\tilde{\pi}-\pi_0)+O_p(n^{-2}). \label{lemma3.01}
\end{equation} 

\noindent {\bf Proof of Lemma 3}\\[0.1in]
It may be verified that under the three conditions and with the composite similarity mapping $h_N^C$ defined by (14) and (22), Theorem 1, Lemma 2 and Theorem 2 all hold. In particular, $\pi_0'-\pi_0=O_p(n^{-3/2})$ and the extended empirical log-likelihood ratio $l^*_2(\pi_0)$ converges in distribution to a $\chi^2_q$ random variable. 

Since $\delta(n)=O(n^{-1/2})$ and $l(\pi_0')=l^*_2(\pi_0)$ which is asymptotically a $\chi^2_q$ variable, we have
\begin{equation}
\{l(\pi_0')\}^{\delta(n)}=1+O_p(n^{-1/2}). \label{lemma3.02}
\end{equation}
By $h^C_N(\pi_0')=\pi_0$, we have $\pi_0-\tilde{\pi}=\gamma_2\{N,l(\pi_0')\}(\pi_0'-\tilde{\pi})$. Thus,
\begin{equation} 
\pi_0'-\pi_0=\frac{b \{l(\pi_0')\}^{\delta(n)}}{2N}(\tilde{\pi}-\pi_0')=\frac{\eta \{l(\pi_0')\}^{\delta(n)}}{2N}(\tilde{\pi}-\pi_0)+\frac{\eta \{l(\pi_0')\}^{\delta(n)}}{2N}(\pi_0-\pi_0'). \label{lemma3.03}
\end{equation}
It follows from (\ref{lemma3.02}), (\ref{lemma3.03}) and  $\pi_0'-\pi_0=O_p(n^{-3/2})$ that
\begin{eqnarray}
\pi_0'-\pi_0&=&\frac{\eta \{l(\pi_0')\}^{\delta(n)}}{2N}(\tilde{\pi}-\pi_0)+O_p(n^{-5/2}) \nonumber \\%\label{3-87} 
                  &=&\frac{\eta}{2N}(\tilde{\pi}-\pi_0)+O_p(n^{-2}), \nonumber %\label{b-87-2} 
\end{eqnarray}
which proves the lemma.
$\Box$ \vspace{0.1in}

% ************************************  Theorem 5

\noindent {\bf Proof of Theorem 5}\\[0.1in]
Under conditions 1, 2, 3 and 4, based on the (\ref{oel.bart}) we can show that $l(\pi_0)$ has the following expansion
\begin{equation}
l(\pi_0)=N(R_1+R_2+R_3)^T(R_1+R_2+R_3) + N\Delta + O_p(n^{-3/2}),\label{t3-01}
\end{equation} 
where $R_i$ and $\Delta$ are functions of $S^{t_1t_2 \dots t_l}$ and $S^{t_1t_2 \dots t_l}_1$ with
\begin{eqnarray}
  R^\tau_1&=&S^\tau, \hspace{0.2in} R^\tau_2=-\frac{1}{2}S^{\tau\upsilon}S^{\upsilon} + \frac{1}{3}s^{\tau\upsilon\omega}S^{\upsilon}S^{\omega}, \label{oel.bart.02}\\ R_3^\tau&=&\frac{3}{8}S^{\tau\upsilon}S^{\upsilon\omega}S^{\omega}-\frac{5}{12}s^{\tau\upsilon\omega}S^{\omega\alpha}S^{\upsilon}S^{\alpha}-\frac{5}{12}s^{\upsilon\omega\alpha}S^{\tau\upsilon}S^{\omega}S^{\alpha} \nonumber\\
  & &+ \frac{4}{9}s^{\tau\upsilon\omega}s^{\omega\alpha\beta}S^{\upsilon}S^{\alpha}S^{\beta}+\frac{1}{3}S^{\tau\upsilon\omega}S^{\upsilon}S^{\omega}-\frac{1}{4}s^{\tau\upsilon\omega\alpha}S^{\upsilon}S^{\omega}S^{\alpha},  \label{oel.bart.03}\\
  \Delta &=& (S^{\tau\upsilon}-S^{\tau\upsilon}_1)S^{\tau}S^{\upsilon}+\frac{2}{3}(S^{\tau\alpha\beta}_1-S^{\tau\alpha\beta})S^{\tau}S^{\alpha}S^{\beta} \label{oel.bart.04},
\end{eqnarray}
where for a vector $P,\hspace{0.06in} P^{r}$ means its $r$th component.  Based on the proofs of (\ref{oel.bart}), we have
\begin{eqnarray}
&(i)&    R_j=O_p(n^{-j/2})\hspace{0.2in} \mbox{for $j=1,2,3$}, \label{oel.bart.05}\\
&(ii)&   D_1= \frac{1}{n}\sum g(y_j,\theta_{y}^0) - \frac{1}{m}\sum g(x_i,\theta_{x}^0) = O_p(n^{-1/2}), \label{oel.bart.06}\\
&(iii)&  \lambda(\pi_0)=V^{-1}D_1+O_p(n^{-1}), \label{oel.bart.07} \\
&(iv)&   R_1^{T}R_1=D_1^{T}V^{-1}D_1, \label{oel.bart.08} \\
&(v)&    \Delta=O_p(n^{-3/2}), \label{oel.bart.09} 
\end{eqnarray}
It may be verified that Lemma 2, Theorem 1 and Theorem 2 all hold under $\gamma_2(N,l(\pi))$. Hence, the limiting distribution of $l^*_2(\pi_0)$ is also $\chi^2_q$. This and the condition that $\delta(n)=O(n^{-1/2})$ imply
\begin{equation}
\left[l(\pi_0')\right]^{\delta(n)}=1+O_p(n^{-1/2}). \label{eel2.02}
\end{equation}  
Since $h_N^C(\pi'_0)=\pi_0$, by (\ref{h.function}) and (\ref{gamma_2}), we have
\begin{eqnarray}
\pi_0'-\pi_0 &=&\frac{\eta \left\{l(\pi_0')\right\}^{\delta(n)}}{2N}(\tilde{\pi}-\pi_0') \nonumber \\ 
                   &=&\frac{\eta \left\{l(\pi_0')\right\}^{\delta(n)}}{2N}(\tilde{\pi}-\pi_0)+\frac{\eta [l(\pi_0')]^{\delta(n)}}{2N}(\pi_0-\pi_0'). \label{eel2.03}
\end{eqnarray} 
By the part ($ii$) of Lemma 2, (\ref{eel2.02}) and (\ref{eel2.03}), we find that
\begin{eqnarray}
\pi_0'-\pi_0 &=& \frac{\eta \left\{l(\pi_0')\right\}^{\delta(n)}}{2N}(\tilde{\pi}-\pi_0)+ O_p(n^{-5/2})\nonumber \\  
                   &=&\frac{\eta}{2N}(\tilde{\pi}-\pi_0)+ O_p(n^{-2}). \label{eel2.04}
\end{eqnarray} 
By (\ref{lemma3.01}) from Lemma 3 and Taylor expansion (\ref{t2-01}), we have 
\begin{eqnarray}
l^*_2(\pi_0)&=& l(\pi_0)+J(\pi_0)(\pi_0'-\pi_0) + o_p(n^{-3/2}) \nonumber \\
              &=& l(\pi_0)+\frac{\eta}{2n}J(\pi_0)(\tilde{\pi}-\pi_0) + O_p(n^{-3/2}), \label{t5-02}
\end{eqnarray} 
where $J(\pi_0)$ is given by (\ref{lemma1.01}). Denote

\begin{equation}
G(X,Y,\pi)=\frac{1}{n}\sum g(y_j,\theta_y) - \frac{1}{m}\sum g(x_i,\theta_x),\label{t5-03}
\end{equation} 
Under condition 2, Taylor expansion of $G(X,Y,\tilde{\pi})$ at $\pi_0$ gives
\begin{eqnarray}
G(X,Y,\tilde{\pi})=&&G(X,Y,\pi_0)+G'(X,Y,\pi_0)(\tilde{\pi}-\pi_0)+O_p(\| \pi_0-\tilde{\pi}\|^2) \nonumber \\
             =&&\left\{\frac{1}{n}\sum g(y_j,\theta_{y_0}) - \frac{1}{m}\sum g(x_i,\theta_{x_0})\right\}+  \nonumber \\
              && \left\{\frac{1}{n}\sum g'(y_j,\theta_{y_0})\frac{\partial \theta_y}{\partial \pi}|_{\pi_0}-\frac{1}{m}\sum g'(x_i,\theta_{x_0})\frac{\partial \theta_x}{\partial \pi}|_{\pi_0}\right\}  \nonumber \\
              &&(\tilde{\pi}-\pi_0)+O_p(n^{-1})  \nonumber 
\end{eqnarray}
This and the estimating equations are just-determined, $G(X,Y,\tilde{\pi})=0$, imply
\begin{eqnarray}
\nonumber && \left\{ \frac{1}{n}\sum g'(y_j,\theta_{y_0})\frac{\partial \theta_y}{\partial \pi}|_{\pi_0}- \frac{1}{m}\sum g'(x_i,\theta_{x_0})\frac{\partial \theta_x}{\partial \pi}|_{\pi_0} \right\} (\pi_0-\tilde{\pi})\nonumber \\
              && =\left\{\frac{1}{n}\sum g(y_j,\theta_{y_0}) - \frac{1}{m}\sum g(x_i,\theta_{x_0})\right\}+  O_p(n^{-1})  \label{t5-04}
\end{eqnarray}
Noting that $\lambda(\pi_0)=O_p(n^{-1/2})$ and  $\pi_0-\tilde{\pi}=O_p(n^{-1/2})$, we can show
\begin{footnotesize}
\begin{eqnarray}
\nonumber &&  \left\{\frac{1}{n} \sum 
              \frac{ g'(y_j,\theta_{y_0})\frac{\partial \theta_y}{\partial \pi}|_{\pi=\pi_0} }{1+f_n \lambda^T(\pi_0)g(y_j,\theta_{y_0})}
              -\frac{1}{m}\sum \frac{g'(x_i,\theta_{x_0})\frac{\partial \theta_x}{\partial \pi}|_{\pi=\pi_0}}{1-f_m \lambda^T(\pi_0)g(x_i,\theta_{x_0})} \right\}  (\pi_0-\tilde{\pi}) \nonumber \\
             && =\left\{ \frac{1}{n}\sum g'(y_j,\theta_{y_0})\frac{\partial \theta_y}{\partial \pi}|_{\pi_0}- \frac{1}{m}\sum g'(x_i,\theta_{x_0})\frac{\partial \theta_x}{\partial \pi}|_{\pi_0} \right\} (\pi_0-\tilde{\pi})+O_p(n^{-1})  \label{t5-05}
\end{eqnarray}
\end{footnotesize}
It follows from (\ref{t5-04}) and (\ref{t5-05}) that
\begin{eqnarray}
\nonumber &&  \left\{\frac{1}{n} \sum 
              \frac{ g'(y_j,\theta_{y_0})\frac{\partial \theta_y}{\partial \pi}|_{\pi=\pi_0} }{1+f_n \lambda^T(\pi_0)g(y_j,\theta_{y_0})}
              -\frac{1}{m}\sum \frac{g'(x_i,\theta_{x_0})\frac{\partial \theta_x}{\partial \pi}|_{\pi=\pi_0}}{1-f_m \lambda^T(\pi_0)g(x_i,\theta_{x_0})} \right\}  (\pi_0-\tilde{\pi}) \nonumber \\
             && = \left\{\frac{1}{n}\sum g(y_j,\theta_{y_0}) - \frac{1}{m}\sum g(x_i,\theta_{x_0})\right\}+  O_p(n^{-1})  \label{t5-06}
 %            && = \left[\frac{1}{n}\sum g(y_j,\theta_{y}^0) - \frac{1}{m}\sum g(x_i,\theta_{x}^0)\right]+  O_p(n^{-1/2})   \label{t5-06}
\end{eqnarray}
By (\ref{t5-02}), (\ref{lemma1.01}) and (\ref{t5-06}), we have
\begin{footnotesize}
\begin{eqnarray}
l^*_2(\pi_0)&=& l(\pi_0)+\frac{\eta}{2N}J(\pi_0)(\tilde{\pi}-\pi_0) + O_p(n^{-3/2}) \nonumber \\
  &=& l(\pi_0)-\frac{\eta}{2N}2 N \lambda^T(\pi_0)\frac{1}{N} \nonumber\\
  &&\left\{ \sum 
              \frac{f_n g'(y_j,\theta_{y_0})\frac{\partial \theta_y}{\partial \pi}|_{\pi=\pi_0} }{1+f_n \lambda^T(\pi_0) g(y_j,\theta_{y_0})}
              -\sum \frac{f_m g'(x_i,\theta_{x_0})\frac{\partial \theta_x}{\partial \pi}|_{\pi=\pi_0}}{1-f_m \lambda^T(\pi_0)g(x_i,\theta_{x_0})} \right\}(\tilde{\pi}-\pi_0)+O_p(n^{-3/2}) \nonumber \\                             
   &=& l(\pi_0)-\frac{\eta}{N}N\lambda^T(\pi_0)\left\{\frac{1}{n}\sum g(y_j,\theta_{y_0}) - \frac{1}{m}\sum g(x_i,\theta_{x_0}) \right\}+O_p(n^{-3/2}) \label{t5-07} 
    % &=&l(\pi_0)-\frac{b}{n}n\lambda^T(\pi_0)
 % \left\{n^{-1}\sum_{i=1}^n g(X_i,\pi_0)+O_p(n^{-1}) \right\}+O_p(n^{-3/2}) \nonumber \\
 % &=&l(\pi_0)-\frac{b}{n}n\lambda^T(\pi_0)
 % \left\{n^{-1}\sum_{i=1}^n g(X_i,\pi_0)\right\}+O_p(n^{-3/2}).\label{b-102}                  
\end{eqnarray}
\end{footnotesize}
Finally, by(\ref{t3-01}), (\ref{t5-07}), and from (\ref{oel.bart.05}) to (\ref{oel.bart.09}), we have
\begin{eqnarray}
l^*_2(\pi_0)&=&l(\pi_0)-\frac{\eta}{N} \left\{N(V^{-1}D_1+O_p(n^{-1}))^T(D_1+O_p(n^{-1}))\right\}+O_p(n^{-3/2})  \nonumber \\
      &=&l(\pi_0)-\frac{\eta}{N} NR_1^TR_1+O_p(n^{-3/2})  \nonumber \\
      &=& l(\pi_0)-\frac{\eta}{N} \{N(R_1+R_2+R_3)^T(R_1+R_2+R_3)+N\Delta \}+O_p(n^{-3/2})  \nonumber \\
      &=& l(\pi_0)-\frac{\eta}{N} l(\pi_0)+O_p(n^{-3/2})  \nonumber \\
      &=& l(\pi_0)\left\{1-\frac{\eta}{N} +O_p(n^{-3/2})\right\},  \nonumber
\end{eqnarray}

which proves Theorem 5.  
$\Box$ \vspace{0.1in}

\noindent {\em Remark.} The second-order result of Theorem 5 holds only for the just-determined case as the proof above used the condition that $G(X,Y,\tilde{\pi})=0$ to obtain (\ref{t5-04}). For over-determined cases, a weaker condition $G(X,Y,\tilde{\pi})=O_p(n^{-1})$ would also allow us to get (\ref{t5-04}). However, it is not clear that outside of the just-determined cases when this weaker condition would hold. When this weaker condition does not hold, the extended empirical log-likelihood ratio $l^*_2(\pi)$ defined in Theorem 5 reduces to a first-order extended empirical log-likelihood ratio as Theorem 4 is still valid for $l^*_2(\pi)$.

\vspace{0.3in}

\noindent {\bf Acknowledgements}

\vspace{0.1in}

%We thank two referees for their valuable comments which have led many improvements in this paper. The work of the second author was supported by the National Science and Engineering Research Council of Canada.

%The data set about the extensibility (\%) at 100gm/cm for different level of quality facbric were given in the article of page 138  
% "\textit{Compatibility of Outer and Fusible Interlining Fabrics in Tailored Garments (Textile Research Journal 1997 67: 137)}". 
%
%\newpage

\vspace{0.3in}

\noindent {\bf References}
%\begin{thebibliography}{7}

\end{document}